\documentclass[10pt, a4paper]{article}
\usepackage{amsmath, amssymb, amsthm, amscd}

\newtheorem {theorem}{Theorem}[section]
\newtheorem{corollary}[theorem]{Corollary}
\newtheorem{lemma}[theorem]{Lemma}
\newtheorem {proposition}[theorem]{Proposition}
\theoremstyle {definition}
\newtheorem {definition}[theorem]{Definition}
\theoremstyle {remark}
\newtheorem {remark}[theorem]{Remark}
\newtheorem {example}[theorem]{Example}

\def\Ann{\operatorname{Ann}}
\def\Ass{\operatorname{Ass}}

\newcommand{\fm}{\ensuremath{\mathfrak m}}
\newcommand{\fp}{\ensuremath{\mathfrak p}}

\newcommand{\un}{\ensuremath{\underline}}
\newcommand{\F}{\ensuremath{\mathcal F}}
\newcommand{\D}{\ensuremath{\mathcal D}}
\newcommand{\ha}{\ensuremath{I_{\mathcal F,M}}}
\newcommand{\scm}{\text{sequentially Cohen-Macaulay }}

\begin{document}
\title{On Sequentially Cohen-Macaulay Modules}
\author{Nguyen Tu Cuong\footnote{Email: ntcuong@math.ac.vn}\ \ and Doan Trung Cuong\footnote{Email: dtcuong@math.ac.vn}\\
Institute of Mathematics\\
18 Hoang Quoc Viet Road, 10307 Hanoi, Vietnam}
\date{ }
\maketitle

\begin{abstract}
In this paper we present  characterizations of sequentially
Cohen-Macaulay modules in terms of systems of parameters, which are
 generalizations of well-known results on Cohen-Macaulay and
generalized Cohen-Macaulay modules.
The sequentially Cohen-Macaulayness of Stanley-Reisner rings of small embedding dimension are also examined.\\
{\it Keywords:} sequentially Cohen-Macaulay module, dimension filtration, good system of parameters,
 Stanley-Reisner ring. \\
{\it AMS Classification:} 13H10, 13C13, 13H15.
\end{abstract}


\section{Introduction}
The concept of sequentially Cohen-Macaulay modules was introduced first by Stanley \cite{st} for graded rings. Similarly one can give  the definition of sequentially Cohen-Macaulay modules on local rings (see \cite {cn}): Let $M$ be a finitely generated module over a local ring $R$ with $d=\dim M$. $M$ is called a sequentially Cohen-Macaulay module, if there exists a filtration of submodules of $M$
$$\mathcal D:\ D_0\subset D_1\subset \cdots \subset D_t=M$$ such that each $D_i/D_{i-1}$ is
Cohen -Macaulay and $$0<\dim (D_{1}/D_0)< \dim (D_2/D_1)<\ldots <\dim (D_t/D_{t-1})=d.$$
 It is clear that an unmixed module $M$ is sequentially Cohen-Macaulay if and only if $M$ is Cohen-Macaulay. Therefore in this case the sequentially Cohen-Macaulayness of $M$ can be characterized by the equality $\ell (M/\un x M)=e(\un x;M)$ for some system of parameters $ \un x=(x_1,\ldots x_d)$ of $M$. Let $t=1$ in the filtration $\D$ above. Then it is easy to show in this case that $M$ is sequentially Cohen-Macaulay if and only if $ H_{\fm}^i(M) =0$ for $i=1,\ldots d-1$, where $H_{\fm}^i(M)$ is the $i^{th}$ local cohomology module of $M$ with respect to the maximal ideal $\fm$. It follows from the theory of generalized Cohen-Macaulay modules (see \cite {t}) that $M$ is sequentially Cohen-Macaulay if and only if there exists a system of parameters $ \un x=(x_1,\ldots x_d)$ of $M$ such that $\ell (M/\un xM)=\ell (D_0) + e(\un x; D_1)$. Therefore it raises to a natural question: Does there exist a characterization of \scm modules in terms of systems of parameters which generalizes that characterizations of special cases as above?

The purpose of this paper is to give an answer for this question. To this end, we present in Section 2 the notions of filtration satisfying the dimension condition and of good system of parameters with
respect to that filtration. Some basic properties of good systems of parameters are given in this section.   The main results of this paper are stated and proved in Section 3. In the last section, we apply the main theorem to study the sequentially Cohen-Macaulayness of Stanley-Reisner rings of small embedding dimension.

\section{Good systems of parameters}
Throughout this paper, $M$ is a finitely generated module over a local Noetherian commutative ring $(R,\fm)$ with $\dim M=d$.
\begin{definition}
(i) We say that a filtration of submodules of $M$
$$\mathcal F : \ M_0\subset M_1\subset \cdots \subset M_t=M,$$
{\it satisfies the dimension condition} if $\dim M_{i-1}<\dim M_i$
for $i=1, \ldots, t$.\\
(ii) A filtration satisfying the dimension condition,
$$\mathcal D:\
D_0\subset D_1\subset \cdots \subset D_t=M,$$
 is called the {\it
dimension filtration} of $M$ if  the following two conditions are
satisfied

a) $D_0=H_\fm^0(M)$ the zero$^{th}$ local cohomology module of $M$
with respect to the maximal ideal $\fm$;

b) $D_{i-1}$ is the largest submodule of $D_i$ with $\dim D_{i-1}<
\dim D_i$ for all $i= t, t-1, \ldots, 1$.
\end{definition}

\begin{definition} Let $\mathcal F : \ M_0\subset M_1\subset \cdots \subset M_t=M$ be
a filtration satisfying the dimension condition and $d_i=\dim
M_i$. A system of parameters  $\un x=(x_1,$ $\ldots, x_d)$ of $M$
is called a {\it good system of parameters} with respect to the
filtration $\mathcal F$ if $M_i\cap (x_{d_i+1},\ldots,x_d)M=0$ for
$i=0, 1, \ldots, t-1$. A good system of parameters with respect to the dimension
filtration is simply called a good system of parameters of $M$.
\end{definition}
\begin{remark}\label{rmq} i) Because of the Noetherian property of $M$, there always exists the dimension
filtration $\mathcal D$ of $M$  and it is unique. Moreover, let
$\bigcap_{\fp\in \text{Ass}(M)}N(\fp)=0$ be a reduced primary
decomposition of the zero module  of $M$. Then
$D_i=\bigcap_{\dim(R/\fp)\geq d_{i+1}}N(\fp)$, where $d_i=\dim
D_i$.\\
ii) Let $N$ be a submodule of $M$ and $\dim N<\dim M$. From the
definition of the dimension filtration, there exists a $D_i$ such
that $N\subseteq D_i$ and $\dim N=\dim D_i$. Therefore, if $\F:\
M_0\subset M_1\subset \cdots \subset M_{t^\prime}=M$ is a
filtration satisfying the dimension condition, for each $M_j$
there exists a $D_i$ such that $M_j\subseteq D_i$ and $\dim
M_j=\dim
D_i$.\\
iii) If a system of parameters $\un x=(x_1,\ldots ,x_d)$ is good with respect to a
filtration $\mathcal F$,  so is $(x_1^{n_1}, \ldots, x_d^{n_d})$
for any positive
integers $n_1,\ldots,n_d$.\\
iv) A good system of parameters of $M$ is also a good system of parameters with respect to any
filtration satisfying the dimension condition.
\end{remark}

The first result of this section is about the existence of good
system of parameters.
\begin{lemma}\label{exist}
There always exists a good system of parameters of $M$.
\end{lemma}
\begin{proof}
Let $\mathcal D : \ D_0\subset D_1\subset \cdots \subset D_t=M$ be
the dimension filtration of $M$ with $d_i=\dim D_i$. By Remark
\ref{rmq}, (i), $D_i=\bigcap_{\dim(R/\fp)\geq d_{i+1}}N(\fp)$
where $\bigcap_{\fp\in \Ass M}N(\fp)=0$ is a reduced primary
decomposition of $0$ of $M$. Put $N_i=\bigcap_{\dim(R/\fp)\leq
d_i}N(\fp).$ Then $D_i\cap N_i=0$ and $\dim (M/N_i)=d_i$. By the
Prime Avoidance Theorem there exists a system of parameters $\un
x=(x_1,\ldots,x_d)$ such that $x_{d_i+1},\ldots,x_d
\in\Ann(M/N_i)$. Therefore $(x_{d_i+1},\ldots,x_d)M\cap
D_i\subseteq N_i\cap D_i=0$ as required.
\end{proof}
The   following result is an immediate consequence of Lemma \ref
{exist} and Remark \ref {rmq}, (ii).
\begin {corollary}
Let $\mathcal F : \ M_0\subset M_1\subset \cdots \subset M_t=M$ be
a filtration satisfying the dimension condition. There always
exists a good system of parameters with respect to $\mathcal F$.
\end {corollary}
 Let $\mathcal F :
\ M_0\subset M_1\subset \cdots \subset M_t=M$ be a filtration of
$M$ satisfying the dimension condition with $d_i=\dim M_i$ and
$\un x=(x_1,\ldots,x_d)$ a good system of parameters with respect to $\mathcal F$.
It is clear that $(x_1,\ldots,x_{d_i})$ is a system of parameters of $M_i$.
Therefore the following difference is well defined
$$\ha(\un x)=\ell(M/\un xM)-\sum_{i=0}^te(x_1,\ldots,x_{d_i};M_i),$$
where $e(x_1,\ldots,x_{d_i};M_i)$ is the Serre multiplicity and we
set $e(x_1,\ldots,x_{d_0};M_0)=\ell(M_0)$ if
 $\dim M_0\leq 0$.
\begin{lemma}\label{positive}
Let $\mathcal F$ be a filtration of $M$ satisfying the dimension
condition and $\un x=(x_1,\ldots,x_d)$ a good system of parameters with respect to
$\mathcal F$. Then $\ha(\un x)\geq 0$.
\end{lemma}
\begin{proof}
Denote
$$\frac{\F}{x_d\F}:\ \frac{M_0+x_dM}{x_dM}\subset\frac{M_1+x_dM}{x_dM}\subset\
 \cdots \subset\frac{M_s+x_dM}{x_dM}\subset\frac{M}{x_dM},$$
where $s=t-1$ if $d_{t-1}<d-1$ and $s=t-2$ if $d_{t-1}=d-1$. Put $\un x^\prime=(x_1,\ldots,x_{d-1})$. Since $(M_i+x_dM)/x_dM\simeq M_i$ for $i\leq s$, the filtration $\mathcal F/x_d\mathcal F$ satisfies the dimension condition and it is easy to prove that $\un x^\prime$ is a good system of parameters of $M/x_dM$ with respect to $\mathcal F/x_d\mathcal F$. On the other hand, we have $$I_{\F/x_d\F,M/x_dM}(\un x^\prime)=\ell(M/\un xM)-e(\un x^\prime;0:_Mx_d)-e(\un x;M)- \sum_{i=0}^se(x_1,\ldots,x_{d_i};M_i).$$
If $d_{t-1}<d-1$, $\ha(\un x)-I_{\F/x_d\F,M/x_dM}(\un x^\prime)=e(\un x^\prime;0:_Mx_d)\geq 0.$
If $d_{t-1}=d-1$, then $M_{t-1}\subseteq 0:_Mx_d$ since
$M_{t-1}\cap x_dM=0$. Hence
$\ha(\un x)-I_{\F/x_d\F,M/x_dM}(\un x^\prime)=e(\un x^\prime;0:_Mx_d)-e(\un x^\prime;M_{t-1})\geq 0.$
Therefore $\ha(\un x)\geq I_{\F/x_d\F,M/x_dM}(\un x^\prime)$ and
the lemma follows immediately by induction on $d$.
\end{proof}

Lemma \ref{positive} leads to some consequences which are useful
in the sequel.
\begin{corollary}\label{lon}
Let $\un x=(x_1,\ldots,x_d)$ be a good system of parameters with respect to the
filtration $\mathcal F$. Then $\ha(\un x)\geq
I_{\F/x_d\F, M/x_dM}(x_1,\ldots,x_{d-1})$.
\end{corollary}
\begin{corollary}\label{duong}
Let $\un x=(x_1,\ldots,x_d)$ be a good system of parameters with respect to the
filtration $\F$. Then
$$e(x_1,\ldots,x_r;M/(x_{r+1},\ldots,x_d)M)\geq \sum_{d_i\geq
r}e(x_1,\ldots,x_{d_i};M_i),$$ for all $r=1,2,\ldots,d$.
\end{corollary}
\begin{proof}
By Lech's formula and Lemma \ref{positive} we have
\[\begin{aligned}
e(x_1,\ldots,x_r;M/(x_{r+1},\ldots,x_d)M)&=\lim_n\frac{1}{n^r}\ell\big(M/(x_1^n,
\ldots,x_r^n,x_{r+1},\ldots,x_d)M\big)\\
&\geq\lim_n\frac{1}{n^r}\sum_{i=0}^te(x_1^n,\ldots,x_r^n,x_{r+1},\ldots,x_{d_i};M_i)\\
&= \sum_{d_i\geq r}e(x_1,\ldots,x_{d_i};M_i)
\end{aligned}\]
as required.
\end{proof}
Set $\un x(\un n)=(x_1^{n_1},\ldots,x_d^{n_d})$ for any $d$-tuple
of positive integers $\un n =(n_1,\ldots,n_d)$. We consider the
difference $\ha(\un x(\un n))$ as a function in $n_1,\ldots,n_d$.

\begin{proposition}
Let $\mathcal F$ be a filtration satisfying the dimension
condition and $\un x=(x_1,\ldots,x_d)$ a good system of parameters with respect to
$\mathcal F$. Then the function $\ha(\un x(\un n))$ is increasing.
\end{proposition}
\begin{proof}
We only need to prove that the function
$\ha(x_1,\ldots,x_{r-1},x_r^n,x_{r+1},\ldots,x_d)$
is increasing in $n$ for each $r\in \{1,2,\ldots,d\}$. Put $\un
x(n)=(x_1,\ldots,x_{r-1},x_r^n,x_{r+1},\ldots,x_d)$. We have
\[\begin{aligned}
\ha(\un x(n+1))-\ha(\un x(n))=&\ell(M/\un x(
n+1)M)-\ell(M/\un x(n)M)\\
&-\sum_{d_i\geq r}e(x_1,\ldots,x_{d_i};M_i)\\
\geq &
e(x_r;M/(x_1,\ldots,x_{r-1},x_{r+1},\ldots,x_d)M)\\&-\sum_{d_i\geq
r}e(x_1,\ldots,x_{d_i};M_i).
\end{aligned}\]
Applying Lech's formula we obtain
\[\begin{aligned}
e(x_r;M/(x_1,\ldots,x_{r-1},x_{r+1},\ldots,x_d)M)
&=\lim_{n} \frac{1}{n}\ell\big(M/(x_1,\ldots,x_{r-1},x_r^n,x_{r+1},\ldots,x_d)M\big)\\
&\geq \lim_{n} \frac{1}{n^r}\ell\big(M/(x_1^n,\ldots,x_r^n,x_{r+1},\ldots,x_d)M\big)\\
&=e\big(x_1,\ldots,x_r;M/(x_{r+1 },\ldots,x_d)M\big).
\end{aligned}\]
Therefore, $\ha(\un x(n+1))\geq \ha(\un x(n))$ by Corollary
\ref{duong}.
\end{proof}

\section{Sequentially Cohen-Macaulay modules}
The aim of this section is to give  characterizations of
sequentially Cohen-Macaulay modules in terms of systems of
parameters.
\begin{definition}
A module $M$ is called a {\it sequentially Cohen-Macaulay module}
if for the dimension filtration $\mathcal D:\ D_0\subset
D_1\subset \cdots \subset D_t=M$, each module $D_i/D_{i-1}$ is
Cohen-Macaulay for $i=1,2,\ldots,t$.
\end{definition}
Note that the notion of sequentially Cohen-Macaulay module was
introduced first by Stanley \cite{st} for graded case (see also
Herzog-Sbara \cite{hs}).
In our study of sequentially Cohen-Macaulay modules, the notion of dd-sequences defined in \cite{c} is used frequently. For convenience, we recall briefly the definition and some basic results of dd-sequences presented in \cite {c}.

Firstly, we recall the definition of d-sequences due to Huneke in \cite{ch}. A sequence $(x_1, x_2,\ldots,x_s)$ of elements of $\fm$ is called a d-{\it sequence} of $M$ if $(x_1,\ldots, x_{i-1})M:~x_j=(x_1,\ldots, x_{i-1})M:x_ix_j$ for $i=1, 2, \ldots, s$ and $j\geq i$.

\begin{definition} \cite[Definition 3.2]{c} A sequence $(x_1, \ldots, x_s)$ is called a dd-{\it sequence} of $M$ if for all positive integers $n_1, \ldots, n_s$ and $i=1, 2, \ldots, s$, the sequence $(x_1^{n_1}, \ldots, x_i^{n_i})$ is a d-sequence of the module $M/(x_{i+1}^{n_{i+1}}, \ldots, x_s^{n_s})M$.
\end{definition}

A dd-sequence has many nice properties, especially when $\un x$ is a system of parameters we have the following characterization.

\begin{lemma}\textup{[4, Corollary 3.6]}\label{char}
Let $\un x$ be a system of parameters of $M$. Then $\un x$ is a dd-sequence of $M$ if and only if there exist integers $a_0, a_1, \ldots, a_d$ such that 
$$\ell\big(M/\un x(\un n)M\big)=\sum_{i=0}^da_in_1\ldots n_i$$
for all positive integers $n_1, \ldots, n_d$. In this case, $$a_i=e\big(x_1, \ldots, x_i; (x_{i+2},\ldots, x_d)M:x_{i+1}/(x_{i+2},\ldots, x_d)M\big).$$
\end{lemma}

In general, it does not require a module to have a system of parameters which is a dd-sequence. By virtue of Lemma \ref{exist}, the  following result proved in \cite{c} shows that a sequentially Cohen-Macaulay module always admits a system of parameters which is a dd-sequence.
\begin{proposition}\textup{(See [4, Theorem 1.5])}\label{seqCM}
Let $M$ be a sequentially Cohen-Macaulay module and $\un x$ a system of parameters of $M$. Then  $\un x$ is a good system of parameters if and only if $I_{\mathcal D,M}(\un x(\un n))=0$ for all positive integers $n_1, \ldots, n_d$, in particular, $\un x$ is a dd-sequence.
\end{proposition}

In order to prove the main theorems we need some auxiliary results. 
\begin{lemma}\label{lemma} Let $\un x=(x_1,\ldots,x_d)$ be a
system of parameters of $M$ and $\mathcal D:\ D_0\subset D_1\subset \cdots \subset
D_t=M$ the dimension filtration of $M$ with $\dim D_i=d_i$.
Suppose that $\un x$ is a dd-sequence of $M$. Then  $D_i=0:_Mx_{d_i+1}$ for
$i=0,1,\ldots,t-1$ and $\un x$ is a good system of parameters of $M$.
\end{lemma}
\begin{proof}
We proved the lemma by induction on the dimension of $M$. The case
$d=1$ is clear. Assume that $d>1$. From the hypothesis $\un
x$ is a dd-sequence, we have $D_{t-1}=0:_Mx_d$ by \cite[Lemma 6.3]{c}, and therefore $D_i\cap x_dM\subseteq D_{t-1}\cap x_dM=0$
for all $i=0,\ldots , t-1$.  Assume now that either $d_{t-1}<d-1$
or $i<t-1$. Then $\dim (D_i+x_d^{n_d}M)/x_d^{n_d}M=d_i$ for all
$n_d\geq 1$, since
$$(D_i+x_d^{n_d}M)/x_d^{n_d}M\simeq D_i/x_d^{n_d}M\cap D_i=D_i.$$ From Remark \ref{rmq},(ii),
there exists an $R$-module $D$ in the dimension filtration of $M/x_d^{n_d}M$ such that
$(D_i+x_d^{n_d}M)/x_d^{n_d}M\subseteq D$ and $\dim D=d_i$. Since
the system of parameters $\un x'=(x_1,\ldots, x_{d-1})$ is also a dd-sequence of the module
$M/x_d^{n_d}M$ and $\dim D=d_i$, it follows from the induction hypothesis that
$D=(0:x_{d_i+1})_{M/x_d^{n_d}M}$. Thus $D_i+x_d^{n_d}M\subseteq
x_d^{n_d}M:x_{d_i+1}$ for all $n_d>0$ and therefore $D_i\subseteq
0:_Mx_{d_i+1}$ by the Krull Intersection Theorem. On the other hand, since $\un x$ is a d-sequence,
$$(x_{d_i+1},\ldots,x_d)(0:_Mx_{d_i+1})=0.$$
Thus $\dim(0:_Mx_{d_i+1})\leq d_i$. This implies by the maximality of $D_i$ that  $D_i=0:_Mx_{d_i+1}$ and the first conclusion of the lemma is prove. To show the second conclusion of
the lemma we note by the induction hypothesis that $\un x'$ is a good system of parameters of $M/x_dM$. Therefore $\un x'$ is a good system of parameters with respect to the filtration
$$\frac{\D}{x_d \D}:\ \frac{D_0+x_dM}{x_dM}\subset\frac{D_1+x_dM}{x_dM}\subset\
 \cdots \subset\frac{D_s+x_dM}{x_dM}\subset\frac{M}{x_dM},$$
 where $s=t-1$ if $d_{t-1}<d-1$ and $s=t-2$ if $d_{t-1}=d-1$. Thus
  $(D_i+x_dM)\cap
(x_{d_i+1},\ldots,x_{d-1},x_d)M=x_dM$ for all $i=0,\ldots , s$.
Keep in mind that $D_{t-1}=0:_Mx_d$ if $d_{t-1}=d-1$.  Hence
$$D_i\cap
(x_{d_i+1},\ldots,x_{d-1},x_d)M\subseteq D_i\cap x_dM=0$$ for all
$i=0,1,\ldots , t-1$ and $\un x$ is a good system of parameters of
$M$.
\end{proof}

\begin{lemma}\label{them} Let $(x_1,\ldots, x_s)$ be a d-sequence of $M$. Then, for $i=1, \ldots, s$, 
$$0:_Mx_i\cap (x_1,\ldots, x_s)M=0:_Mx_i\cap (x_1,\ldots, x_{i-1})M.$$
In particular, $0:_Mx_2\cap (x_1,\ldots, x_s)M=x_1(0:_Mx_2)$.
\end{lemma}
\begin{proof}
For the case $i=s$, we need only to prove that 
$$0:_Mx_s\cap (x_1,\ldots, x_s)M\subseteq 0:_Mx_s\cap (x_1,\ldots, x_{s-1})M.$$ 
Let $a\in 0:_Mx_s\cap (x_1,\ldots, x_s)M$, then $ax_s=0$ and $a=x_1a_1+\cdots+ x_sa_s$ where $a_1, \ldots, a_s\in M$. Thus $a_s\in (x_1,\ldots, x_{s-1})M:x_s^2=(x_1,\ldots, x_{s-1})M:x_s$ and so $a \in 0:_Mx_s\cap (x_1,\ldots, x_{s-1})M$.
The case $i<s$ is proved by decreasing induction. Assume that $0:_Mx_{i+1}\cap (x_1,\ldots, x_s)M=0:_Mx_{i+1}\cap (x_1,\ldots, x_i)M$. Since $0:_Mx_i\subseteq 0:_Mx_{i+1}$, we have
$0:_Mx_i\cap (x_1,\ldots, x_s)M=0:_Mx_i\cap (x_1,\ldots, x_i)M$. Note that $(x_1,\ldots, x_i)$ is also a d-sequence, then by our proof above for $i=s$, we have $0:_Mx_i\cap (x_1,\ldots, x_s)M=0:_Mx_i\cap (x_1,\ldots, x_{i-1})M$. 
\end{proof}
The following result is a key lemma for the proof of our main results
.
\begin{lemma}\label{key}
Let $\mathcal D:\ D_0\subset D_1\subset \cdots \subset D_t=M$ be
the dimension filtration of $M$ and $\mathcal F:\ M_0\subset
M_1\subset \cdots \subset M_{t^\prime}=M$ a filtration satisfying
the dimension condition. Assume that $\un x=(x_1,\ldots,x_d)$ is a
good system of parameters with respect to $\mathcal F$ and $\ha(\un x(\un n))=0$
for all positive integers $n_1,\ldots,n_d$. Then $I_{\mathcal
D,M}(\un x(\un n))=0$, $t=t^\prime$ and $\dim M_i=\dim D_i=d_i$.
Moreover, we have $D_i+x_1M=x_1M:x_{d_i+1}$ for $i=1,\ldots ,
t-1$.
\end{lemma}
\begin{proof}
Since $\ha(\un x(\un n))=0$, $\un x$ is dd-sequence (Lemma \ref{char}), so is a good system of parameters of $M$ by Lemma \ref{lemma}. From Remark \ref{rmq},(ii), we have $\ha(\un x(\un n))\geq I_{\mathcal D,M}(\un x(\un n))$ for all positive integers $n_1,\ldots,n_d$. Therefore $I_{\mathcal D,M}(\un x(\un n))=0$, $t=t^\prime$ and
$\dim D_i=\dim M_i$. We proceed  the last conclusion by induction on $d$. The case $d=1$ is trivial. Let $d \geq 2$. Assume first that $\dim D_1>1$. Then $x_2$ is a parameter element of $D_i$ for
$i=1,\ldots , t$. For any positive integer $n_2$ we consider the following filtration
$$\frac{\mathcal D}{x_2^{n_2}\mathcal D}:\ \frac{D_0+x_2^{n_2}M}{x_2^{n_2}M}
\subset\frac{D_1+x_2^{n_2}M}{x_2^{n_2}M}\subset\ \cdots \subset\frac{D_{t-1}+x_2^{n_2}M}{x_2^{n_2}M}\subset\frac{M}{x_2^{n_2}M}.$$
By Lemma \ref{lemma}, we have $D_i=0:_Mx_{d_i+1}$ for all $i=1,\ldots,t$,  then $D_i\cap x_2^{n_2}M=x_2^{n_2}D_i$ and 
$\dim\big((D_i+x_2^{n_2}M)/x_2^{n_2}M\big)=\dim(D_i/x_2^{n_2}D_i)=d_i-1$.
Thus the filtration $\mathcal D/x_2^{n_2}\mathcal D$ satisfies the dimension condition. Moreover, it can be proved similarly as Lemma \ref{positive} that $I_{\mathcal D,M}(\un x(\un n))\geq I_{\mathcal D/x_2^{n_2}\mathcal D,M/x_2^{n_2}M}(x_1^{n_1},x_3^{n_3},\ldots,x_d^{n_d})$, which implies $I_{\mathcal D/x_2^{n_2}\mathcal D,M/x_2^{n_2}M}(x_1^{n_1},x_3^{n_3},\ldots,x_d^{n_d})=0$ for all positive integers $n_1,\ldots,n_d$. Let
$\mathcal D' : D'_0\subset D'_1\subset \ldots \subset D'_{t'}=M/x_2^{n_2}M$
be the dimension filtration of $M/x_2^{n_2}M$. Then $I_{\mathcal
D^\prime,M/x_2^{n_2}M}(x_1^{n_1},x_3^{n_3},\ldots,x_d^{n_d})=0$
and $t'=t$ by the first conclusion of the lemma. Therefore,  from Lemma \ref{lemma},
$D'_i=(0:x_{d_i+1})_{M/x_2^{n_2}M}$. Hence
$$(0:x_{d_i+1})_{M/x_2^{n_2}M}+x_1(M/x_2^{n_2}M)=x_1(M/x_2^{n_2}M):x_{d_i+1}$$
by the induction hypothesis; so
$x_2^{n_2}M:x_{d_i+1}+x_1M=(x_1,x_2^{n_2})M:x_{d_i+1}$ for
$i=1,\ldots , t-1$. Applying Lemma \ref{lemma} and the Krull Intersection Theorem, we get
$D_i+x_1M=x_1M:x_{d_i+1}$ for $i=1,\ldots , t-1$.

Now assume that $\dim D_1=1$. Using Corollary 4.3 of
Auslander-Buchsbaum \cite{ab} to the sequence $x_2^{n_2}, x_1^{n_1}$
with respect to the module $M/(x_3^{n_3},\ldots,x_d^{n_d})M$ we have
\[\begin{aligned} 
\ell(M/\un x(\un n))=&e(x_2^{n_2},x_1^{n_1};M/(x_3^{n_3},\ldots,x_d^{n_d})M)
+n_2e(x_2;(0:x_1^{n_1})_{M/(x_3^{n_3},\ldots,x_d^{n_d})M})\\
&+\ell((0:x_2^{n_2})_{M/(x_1^{n_1},x_3^{n_3},\ldots,x_d^{n_d})M}).
\end{aligned}\]
On the other hand, it implies by Lech's formula and the condition $I_{\mathcal
D,M}(\un x(\un n))=0$ that
\[\begin{aligned}
e(x_2^{n_2},x_1^{n_1};M/(x_3^{n_3},\ldots,x_d^{n_d})M)&=\lim_{r}\frac{1}{r^2}
\ell(M/(x_1^{rn_1},x_2^{rn_2},x_3^{n_3},\ldots,x_d^{n_d})M)\\
&=\sum_{i\geq 2}n_1\ldots n_{d_i}e(x_1,\ldots,x_{d_i};D_i)\\
&=\ell(M/\un x(\un n)M)-n_1e(x_1;D_1)-\ell(D_0).
\end{aligned}\]
Therefore,
$$n_1e(x_1;D_1)+\ell(D_0)=n_2e(x_2;(0:x_1^{n_1})_{M/(x_3^{n_3},\ldots,x_d^{n_d})M})
+\ell((0:x_2^{n_2})_{M/(x_1^{n_1},x_3^{n_3},\ldots,x_d^{n_d})M}).$$
Since the left term of the  equality above is independent of $n_2$
and the right term is increasing in $n_2$, we deduce that
$e(x_2;(0:x_1^{n_1})_{M/(x_3^{n_3},\ldots,x_d^{n_d})M})=0$ and
$$\ell (D_1/x_1D_1)=n_1e(x_1;D_1)+\ell(D_0)
=\ell((0:x_2)_{M/(x_1^{n_1},x_3^{n_3},\ldots,x_d^{n_d})M}).$$
Keep in mind that $D_1=0:_Mx_2$ and the sequence $\un x$ is a d-sequence. We obtain by Lemma~\ref {them},
$$(x_1,x_3^{n_3},\ldots,x_d^{n_d})M\cap D_1=x_1M\cap D_1=x_1D_1,$$
thus $D_1/x_1D_1\simeq (D_1+(x_1,x_3^{n_3},\ldots,x_d^{n_d})M)/(x_1,x_3^{n_3},\ldots,x_d^{n_d})M$.
This  implies 
$$D_1+(x_1,x_3^{n_3},\ldots,x_d^{n_d})M=(x_1,x_3^{n_3},\ldots,x_d^{n_d})M:x_2.$$
Applying again the Krull Intersection Theorem we get
$D_1+x_1M=x_1M:x_2$. Now we set $M^\prime=M/D_1$. Then the
dimension filtration of $M^\prime$ is just
$$0\subset D_2/D_1\subset \cdots \subset D_t/D_1=M^\prime.$$
By Lemma \ref{them}, we have $D_1\cap \un x(\un
n)M=x_1^{n_1}D_1$. Hence there is a short exact
sequence
$$0\longrightarrow D_1/x_1^{n_1}D_1\longrightarrow M/\un x(\un n)M
\longrightarrow M^\prime/\un x(\un n)M^\prime\longrightarrow 0.$$
Therefore
$$\ell(M^\prime/\un x(\un n)M^\prime)=\sum_{i=2}^tn_1\ldots
n_{d_i}e(x_1,\ldots,x_{d_i};D_i/D_1).$$ So that the sequence
$(x_1,\ldots,x_d)$ is  satisfied the hypothesis of the lemma with
respect to the module $M^\prime$. Since $\dim(D_2/D_1)>1$, we have
$D_i/D_1+x_1(M/D_1)=x_1(M/D_1):x_{d_i+1}$ as in the proof above.
Therefore $D_i+x_1M=x_1M:x_{d_i+1}$ for all $i\geq 2$.
\end{proof}

We  go now to the first characterization of sequentially Cohen-Macaulay modules by dd-sequence and  the vanishing of function $I_{\mathcal D, M}(\un x)$.
\begin{theorem}\label{dd}
Let $M$ be a finitely generated $R$-module of dimension $d$ and $\mathcal D\ :\ D_0\subset D_1\subset \cdots \subset D_t=M$ the dimension filtration of $M$. Denote $\dim D_i=d_i$. Then the following statements are equivalent.
\begin{itemize}
\item[i)] $M$ is a sequentially Cohen-Macaulay module.

\item[ii)] $I_{\mathcal D,M}(\un x)=0$ for all systems of parameters $\un x$ of $M$ which are dd-sequences.

\item[iii)] There exists a system of parameters $\un x$ of $M$ such that $\un x$ is a dd-sequence and $I_{\mathcal D, M}(\un x)=0$.
\end{itemize}
\end{theorem}
\begin{proof}
$(i \Rightarrow ii)$ follows from Proposition \ref{seqCM} and $(ii \Rightarrow iii)$ is trivial.

\noindent $(iii \Rightarrow i)$.  Let $\un x$ be a system of parameters of $M$ such that $\un x$ is a dd-sequence and $I_{\mathcal D,M}(\un x)=0$. We show first that $I_{\mathcal D, M}(\un x(\un n))=0$ for all positive integers $n_1,\ldots, n_d$. In fact, since  $\un x$ is a good system of parameters of $M$ and $D_i=0:_Mx_{d_i+1} $ by Lemma \ref{lemma}, we have $$D_i\simeq (D_i+ (x_{d_i+2},\ldots, x_d)M)/(x_{d_i+2},\ldots, x_d)M\subseteq (x_{d_i+2},\ldots, x_d)M:x_{d_i+1}\big/(x_{d_i+2},\ldots, x_d)M.$$
Then $e(x_1, \ldots, x_{d_i};D_i)\leq e(x_1,\ldots, x_{d_i};(x_{d_i+2},\ldots, x_d)M:x_{d_i+1}\big/(x_{d_i+2},\ldots, x_d)M)$. Since $\un x$ is a dd-sequence, by Lemma \ref{char} and the hypothesis we have
$$\ell(M/\un xM)=\sum_{i=0}^de(x_1,\ldots, x_i; (x_{i+2},\ldots,x_d)M:x_{i+1}\big/(x_{i+2},\ldots, x_d)M) = \sum_{i=0}^te(x_1,\ldots, x_{d_i}; D_i).$$
Therefore
 $$e(x_1,\ldots, x_i;(x_{i+2},\ldots, x_d)M:x_{i+1}\big/(x_{i+2},\ldots, x_d)M)=0$$ for all $i\notin \{d_0, d_1, \ldots, d_t\}$ and 
$$e(x_1, \ldots, x_{d_i};D_i)=e(x_1,\ldots, x_{d_i};(x_{d_i+2},\ldots, x_d)M:x_{d_i+1}\big/(x_{d_i+2},\ldots, x_d)M).$$
Thus  $I_{\mathcal D, M}(\un x(\un n))=0$ for all positive integers $n_1,\ldots, n_d$ by Lemma \ref{char}. Next, we prove  by induction on $d$ that $M$ is sequentially Cohen-Macaulay. It is clear when $d=1$. By  passing $M$ to the module $M/D_0$ we may assume  that $\dim D_0=0$. If $\dim D_1=1$, $D_1$ is Cohen-Macaulay; therefore $M$ is \scm if and only if so is $M/D_1$. Similarly as in the proof of Lemma \ref{key}, by passing $M$ to $M/D_1$ we may assume in addition
that $\dim D_1>1$.  Consider the following filtration
$$\mathcal D/x_1\mathcal D:\ 0\subset (D_1+x_1M)/x_1M\subset \cdots
 \subset (D_{t-1}+x_1M)/x_1M\subset M/x_1M.$$
We have $D_i\cap x_1M=x_1D_i$ and $\dim((D_i+x_1M)/x_1M)=\dim(D_i/x_1D_i)=d_i-1$. Hence, the
filtration $\mathcal D/x_1\mathcal D$ satisfies the dimension condition. On the other hand, since $D_i+x_1M=x_1M:x_{d_i+1}$ by Lemma \ref{key}, we have $(D_i+x_1M)/x_1M=(0:x_{d_i+1})_{M/x_1M}$ for all $i=1,\ldots , t-1$.  Therefore,
$$I_{\D/x_1\D,M/x_1M}(x_2^{n_2},\ldots,x_d^{n_d})=I_{\D,M}(x_1, x_2^{n_2},\ldots,x_d^{n_d})=0,$$
and the filtration  $\mathcal
D/x_1\mathcal D$ is the dimension filtration of $M/x_1M$ by
Lemmas \ref{lemma}, \ref{key}. Thus  $M/x_1M$ is sequentially
Cohen-Macaulay by the induction hypothesis. It follows that
$$(D_i+x_1M)/(D_{i-1}+x_1M)\simeq (D_i/D_{i-1})/x_1(D_i/D_{i-1})$$
 is Cohen-Macaulay for $i=1,\ldots,t$. So $D_i/D_{i-1}$ is Cohen-Macaulay for $i=1,\ldots,t,$ since
$x_1$ is a regular element of $D_i/D_{i-1}$. Therefore $M$ is a sequentially Cohen-Macaulay module.
\end{proof}

The next theorem is a characterization of sequentially Cohen-Macaulay modules in terms of good systems of parameters.
\begin{theorem}\label{main}
The following conditions  are equivalent:
\begin{itemize}
\item[i)] $M$ is a sequentially Cohen-Macaulay module.

\item[ii)] There exists a filtration satisfying the dimension condition $\F$ such that $\ha(x_1^2,\ldots,x_d^2)=~0$ for all good systems of parameters $\un x=(x_1,\ldots , x_d)$ with respect
to $\mathcal F$.

\item[iii)] There exists a filtration satisfying the dimension condition $\F$  such that $\ha(\un x(\un n))=0$ for all good systems of parameters $\un x=(x_1,\ldots , x_d)$ with respect to $\mathcal F$ and all positive integers $ n_1,\ldots, n_d$.

\item[iv)] There exists a filtration satisfying the dimension condition $\F$ and a good system of parameters $\un x=(x_1,\ldots, x_d)$ with respect to $\mathcal F$ such that 
$\ha(x_1^2,\ldots,x_d^2)=0.$

\item[v)] There exists a filtration satisfying the dimension
condition $\F$ and a good system of parameters $\un x=(x_1,\ldots
, x_d)$ with respect to $\mathcal F$ such that $\ha(\un x(\un
n))=0$ for all positive integers $ n_1,\ldots, n_d$.
\end{itemize}
\end{theorem}
\begin{proof}
($i\Rightarrow iii$) follows from Proposition \ref{seqCM}. The implications ($iii\Rightarrow ii$) and ($ii\Rightarrow iv$) are trivial.

\noindent ($v \Rightarrow i$).  Since $\un x$ is a dd-sequence by the hypothesis and Lemma \ref{char}, $\un x$ is a good system of parameters of $M$ by Lemma \ref {lemma}. Then the implication follows  from Proposition \ref{seqCM} and Theorem \ref{dd}.

\noindent ($iv \Rightarrow v$).   Let $\un x=(x_1,\ldots,x_d)$ be a
good system of parameters of $M$ with respect to $\F$ such that
$\ha(x_1^2,\ldots,x_d^2)=0$. Since the function $\ha(\un x(\un
n))$ is increasing and non-negative,  $\ha(\un x(\un n))=0$ for
all $ n_1,\ldots, n_d\in \{  1,2\}$. First we prove that $\ha(\un
x(\un n))=0$ for all $ n_1,\ldots, n_{d-1}\in \{  1,2\}$ and
arbitrary positive integer $n_d$. In fact, applying Corollary 4.3
of  \cite{ab} to the sequence
$x_d^{n_d},x_1^{n_1},\ldots,x_{d-1}^{n_{d-1}}$ we have
\[\begin{aligned}
\ell(M/\un x(\un n)M)-e(\un x(\un
n);M)=&n_d\sum_{i=0}^{d-2}e(x_d,x_1^{n_1},
\ldots,x_i^{n_i};(0:x_{i+1}^{n_{i+1}})_{M/(x_{i+2}^{n_{i+2}},\ldots,x_{d-1}^{n_{d-1}})M})\\
&+\ell((0:x_d^{n_d})_{M/(x_1^{n_1},\ldots,x_{d-1}^{n_{d-1}})M})\\
=&\sum_{i=0}^{t-1}n_1\ldots n_{d_i}e(x_1,\ldots, x_{d_i};M_i),
\end{aligned}\]
for all $n_1,\ldots,n_d\in \{1,2\}$. Since the right term of the
last equality is independent of $n_d$, it follows that
$$\sum_{i=0}^{d-2}e(x_d,x_1^{n_1},\ldots,x_i^{n_i};(0:x_{i+1}^{n_{i+1}})_{M/(x_{i+2}^{n_{i+2}},
\ldots,x_{d-1}^{n_{d-1}})M})=0$$
and
$$(0:x_d^{2})_{M/(x_1^{n_1},\ldots,x_{d-1}^{n_{d-1}})M}=(0:x_d)_{M/(x_1^{n_1},
\ldots, x_{d-1}^{n_{d-1}})M}.$$ Thus
\[\begin{aligned}
\ell(M/\un x(\un n)M)-e(\un x(\un n);M)&=\ell((0:x_d^{n_d})_{M/(x_1^{n_1},\ldots,x_{d-1}^{n_{d-1}})M})
&=\ell((0:x_d)_{M/(x_1^{n_1},\ldots,x_{d-1}^{n_{d-1}})M}),
\end{aligned}\]
is independent of $n_d$. Therefore $\ha(\un x(\un n))=0$ for all
$n_1,\ldots,n_{d-1}\in\{1,2\}$ and all $n_d\geq 1$. Next, we prove
the statement by induction on the dimension of $M$. The case $d=1$
is done by the proof above. Assume $d>1$. For an arbitrary
positive integer $n_d$ we denote
$$\frac{\F}{x_d^{n_d}\F}:\
\frac{M_0+x_d^{n_d}M}{x_d^{n_d}M}\subset\frac{M_1+x_d^{n_d}M}{x_d^{n_d}M}\subset\
\cdots
\subset\frac{M_s+x_d^{n_d}M}{x_d^{n_d}M}\subset\frac{M}{x_d^{n_d}M},$$
 where $s=t-1$ if $d_{t-1}<d-1$ and $s=t-2$ if $d_{t-1}=d-1$. It follows from
Corollary \ref{lon} and the hypothesis that
$I_{\F/x_d^{n_d}\F,M/x_d^{n_d}M}(x_1^{n_1},\ldots,x_{d-1}^{n_{d-1}})=0$
for all $n_1, \ldots, n_{d-1}\in\{1,2\}$. Hence, by the induction
hypothesis,
$I_{\F/x_d^{n_d}\F,M/x_d^{n_d}M}(x_1^{n_1},\ldots,x_{d-1}^{n_{d-1}})=
0,$ for all positive integers $n_1,\ldots,n_{d-1}$. Therefore
\[\begin{aligned}
\ell(M/\un x(\un n)M)=&n_1\ldots
n_{d-1}e(x_1,\ldots,x_{d-1};M/x_d^{n_d}M)+\sum_{i=0}^sn_1
\ldots n_{d_i}e(x_1,\ldots,x_{d_i};M_i)\\
=&n_1\ldots n_de(\un x;M)+n_1\ldots n_{d-1}e(x_1,\ldots,x_{d-1};0:_Mx_d^{n_d})\\
&+\sum_{i=0}^sn_1\ldots
n_{d_i}e(x_1,\ldots,x_{d_i};M_i)\\
\geq & \sum_{i=0}^tn_1\ldots n_{d_i}e(x_1,\ldots,x_{d_i};M_i)
\end{aligned}\]
by Lemma \ref{positive}. Note that the last inequality becomes
an equality for all $n_1,\ldots,n_{d-1}\in \{1,2\}$ and all
positive integers $n_d$.  This implies that
$e(x_1,\ldots,x_{d-1};0:_Mx_d^{n_d})=0$ if $s=t-1$ and
$e(x_1,\ldots,x_{d-1};0:_Mx_d^{n_d})=e(x_1,\ldots,x_{d-1};M_{t-1})$
if $s=t-2$. Thus, it is an equality for all positive integers
$n_1,\ldots,n_d$ and $\ha(\un x(\un n))=0$ as required.
\end{proof}

In \cite{g}, Goto had introduced the notion of approximately Cohen-Macaulay rings. A local ring $(R,\fm)$ is called approximately Cohen-Macaulay if $R$ is not a Cohen-Macaulay ring and there exists an element $a\in \fm$ such that $R/a^nR$ is  a Cohen-Macaulay ring of dimension $d-1$ for every $n>0$. Similarly,  we can define the notion of  an approximately Cohen-Macaulay module.
\begin{definition}
A non Cohen-Macaulay module $M$ is called an {\it approximately Cohen-Macaulay module} if there exists an element $a\in \fm$ such that $M/a^nM$ is Cohen-Macaulay of dimension $d-1$ for every $n>0$.
\end{definition}
 
Then, the following characterization of approximately Cohen-Macaulay modules are 
easily derived by Theorem \ref{main}. It should be mentioned that the equivalence of (i) and (ii) was proof in  \cite [Theorem 1]{g}  for local rings.
\begin{proposition}\label{app}
Let $M$ be not a Cohen-Macaulay $R$-module of dimension $d$. The following statements are equivalent.
\begin{itemize}
\item[i)] $M$ is an approximately Cohen-Macaulay module.

\item[ii)] There exists an element $a\in \fm$ such that $0:_Ma=0:_Ma^2$ and $M/a^2M$ is a Cohen-Macaulay module of dimension $d-1$.

\item[iii)] 
 $M$ is a sequentially Cohen-Macaulay module with the dimension filtration $\mathcal D\ :\ 0=D_0\subset D_1\subset D_2=M$ where $\dim D_1=d-1$.
\end{itemize}
\end{proposition}
\begin{proof}
$(i\Rightarrow ii)$ is trivial.

\noindent $(ii\Rightarrow iii)$. 
Assume that $M/a^2M$ is Cohen-Macaulay of dimension $d-1$. Then there exists a system of parameters $\un x=(x_1, \ldots, x_{d-1}, x_d)$ of $M$ such that $x_d=a$ and we have
\[\begin{aligned}
\ell(M/(x_1^2,\ldots,x_d^2)M)&=e(x_1^2,\ldots,x_{d-1}^2;M/x_d^2M)\\
&=2^de(\un x;M)+2^{d-1}e(x_1,\ldots, x_{d-1}; 0:_Mx_d).
\end{aligned}\]
Since $M$ is not Cohen-Macaulay,
$ e(x_1,\ldots, x_{d-1}; 0:_Mx_d)>0$, so  $\dim (0:_Mx_d)=d-1$. Therefore the filtration $\mathcal F\ :\ 0\subset 0:_Mx_d\subset M$  satisfies the dimension condition. Since $0:_Mx_d=0:_Mx_d^2$, it is easy to prove that $\un x$ is a good system of parameters with respect to the filtration $\mathcal F$. Moreover, we have $\ha(x_1^2, \ldots, x_d^2)=0$. Then $M$ is a sequentially Cohen-Macaulay module by Theorem \ref{main} and therefore $\mathcal F$  is just the dimension filtration of $M$ by Lemma \ref {key}.

\noindent $(iii\Rightarrow i)$. Assume that $M$ is a sequentially Cohen-Macaulay module with the dimension filtration $\mathcal D\ : \ 0=D_0\subset D_1\subset D_2=M$ where $\dim D_1=d-1$. Let $\un x$ be a good system of parameters of $M$. By Proposition \ref{seqCM}, $\un x$ is a dd-sequence and $I_{\mathcal D, M}(\un x(\un n))=0$ for all positive integers $n_1,\ldots, n_d$.  Hence $\ell(M/\un x(\un n)M)=e(x_1^{n_1},\ldots, x_{d-1}^{n_{d-1}};M/x_d^{n_d}M)$, since $D_1=0:_Mx_d=0:_Mx_d^{n_d}$. Therefore $M/x_d^{n_d}M$ is a Cohen-Macaulay module of dimension $d-1$ for all $n_d>0$.
\end{proof}
\begin {remark}
A filtration satisfying the dimension condition
$$\mathcal F : \ M_0\subset M_1\subset \cdots \subset M_t=M$$
 is called a Cohen-Macaulay filtration if $M_i/M_{i-1}$ is Cohen-Macaulay modules for all $i=1,\ldots , t$. Then it was showed in \cite{cn} that if $M$ admits a Cohen-Macaulay filtration $\F$, $M$ is a sequentially Cohen-Macaulay module and $\F$ is just the dimension filtration of $M$. Here we want to clarify that there exists filtration satisfying the equivalent conditions of Theorem \ref{main} which is not the dimension filtration of $M$. Let $M$ be a \scm
module of positive depth and
$$\mathcal D : 0= D_0\subset D_1\subset \cdots \subset D_t=M$$
 the dimension filtration of $M$. Let $\un x=(x_1,\ldots , x_d)$ be a good system of parameters of $M$. Then it is not difficult to check that the following filtration
$$\mathcal F\ :\  0=M_0\subset  M_1=\un x D_1\subset \cdots \subset M_{t-1}=\un x D_{t-1} \subset
M$$ satisfies the equivalent conditions of Theorem 3.4 but it is not the dimension filtration of $M$.
\end{remark}

 It should be mentioned that an $R$-module $M$ is Cohen-Macaulay if and only if there exists a system of parameters (and therefore for all systems of parameters) $\un x$ of $M$ such that $\ell(M/\un xM)-e(\un x;M)=0$. Then it raises to a natural question: a  module  $M$ is sequentially Cohen-Macaulay if and only if there exists a good system of parameters $\un x$ such that $I_{\mathcal D, M}(\un x)=0$ where $\mathcal D$ is the dimension filtration of $M$?  Unfortunately, the answer is negative. Below we give two counter-examples for this question. The first example is due to Goto \cite[Remark 2.9]{g}.
\begin{example} (1)
Let $R=k[[x,y,z,w]]$ be the ring of formal power series over a field $k$ and $P=(xw-yz, x^3-z^2, w^2-xy^2, zw-x^2y)$, $Q=(y^2, z, w)$. Put $M=R/P\cap Q$. Then the dimension filtration of $M$ is $\mathcal D\ :\ 0=D_0\subset D_1\subset D_2=M$ where $D_1=P/P\cap Q$, $\dim D_1=1$. Since $D_1=0:_Mw=0:_Mw^2$, it is easy to see that $(x+y+z+w, w)$ is a good system of parameters of $M$ . Moreover, by a simple computation we have
\begin{equation*}
I_{\mathcal D, M}((x+y+z+w)^{n_1}, w^{n_2})=
\begin{cases}
0&\text{  if   } n_1=n_2=1\\
1&\text{   otherwise}.
\end{cases}
\end{equation*}
Thus  though $I_{\mathcal D, M}(x+y+z+w, w)=0$, $M$ is not a sequentially Cohen-Macaulay module by Theorem \ref {main}, (iii).

\noindent (2) Let $R=k[[x,y,z,w]]$ be the ring of formal power series over a field $k$ and $P=(x,w)\cap(y,z)$, $Q=(x, y^2, z)$. Put $M=R/P\cap Q$. Then the dimension filtration of $M$ is $\mathcal D\ :\ 0=D_0\subset D_1\subset D_2=M$ where $D_1=P/P\cap Q$, $\dim D_1=1$. Since $D_1=0:_M(x+y)=0:_M(x+y)^2$, it is easy to see  that $(z+w, x+y)$ is a good system of parameters of $M$. Then it is easy to check that
\begin{equation*}
I_{\mathcal D, M}((z+w)^{n_1}, (x+y)^{n_2})=
\begin{cases}
0&\text{  if   } n_2=1\\
1&\text{  if   } n_2\geq 2.
\end{cases}
\end{equation*}
So though $I_{\mathcal D, M}(z+w, x+y)=0$, $M$ is not a sequentially Cohen-Macaulay module by Theorem \ref {main}, (iii).
\end{example}

\section{Stanley-Reisner rings}
A simplicial complex $\Delta$ over $n$ vertices
$\{v_1,\ldots,v_n\}$ is a collection of subsets of the set
$\{v_1,\ldots,v_n\}$ such that,

i) $\emptyset \in \Delta$.

ii) For all element $F\in \Delta$ and all subsets $F^\prime\subseteq F$, $F^\prime\in \Delta$.\\
Each element $F\in \Delta$ is called a face of $\Delta$. Among the
faces of $\Delta$, the face $F$ with the property that if
$F\subseteq F^\prime$ and $F^\prime\in \Delta$ then $F=F^\prime$
is called a facet of $\Delta$. So a simplicial complex is defined
completely if all its facets are given.  For a set of $n$ vertices
$\{v_1,\ldots,v_n\}$ we consider the polynomial ring
$R=k[X_1,\ldots,X_n]$ over a field $k$. Then $\Delta$ corresponds
to an ideal $I_\Delta$ of $R$ defined by a set of generators
$\{X_{i_1}\ldots X_{i_s}: \{v_{i_1},\ldots ,v_{i_s}\}\notin
\Delta\}$. The Stanley-Reisner ring of $\Delta$ over $k$ is
defined by $k[\Delta]=R/I_\Delta$. For each face $F$ of $\Delta$
we define $\dim F=\dim R/I_F-1$ and the corresponding ideal $I_F$
is a prime ideal generated by $X_i$'s.  It is well-known that
$I_\Delta$ is a radical ideal and $I_\Delta=\bigcap_{F}I_F$, where
$F$ runs through over all the facets of $\Delta$. Therefore, the
dimension filtration
$$\mathcal D: \ D_0\subset D_1\subset
\cdots \subset D_{t-1}\subset D_t=k[\Delta]$$
 of $k[\Delta]$, $d_i=\dim D_i$,  can
be determined by $D_i=\cap_{\dim F+1\geq d_{i+1}}I_F/I_\Delta.$

For more details on Stanley-Reisner rings, the readers can find in
books of Bruns-Herzog \cite{bh} and Stanley \cite{st}.

Denote by $\lambda_i$ the number of the facets of dimension
$d_i-1$ of the simplicial complex $\Delta$. We derive from Theorem
3.4 the following criterion,  which is very convenient for
checking whether $k[\Delta]$ is a sequentially Cohen-Macaulay
ring.
\begin{proposition}\label{face}
Let $\Delta$ be a simplicial complex on $n$ vertices. Then the
Stanley-Reisner ring $k[\Delta]$  is sequentially Cohen-Macaulay
if and only if there exists a homogeneous good system of
parameters $\un x$ such that
$$\ell(k[\Delta]/(x_1^2,\ldots,x_d^2)k[\Delta])=\sum_{i=1}^t2^{d_i}\lambda_i \deg(x_1)\ldots \deg(x_{d_i}),$$
where $\deg(x_i)$ is the degree of  $x_i$ in the graded ring
$k[\Delta]$.
\begin{proof}
Let $\un x$ be a homogeneous good system of parameters of $M$.
Then the proposition follows from Theorem \ref{main} if  we can
show that
$$e(x_1,\ldots,x_{d_i};D_i)=\lambda_i \deg(x_1)\ldots \deg(x_{d_i}).$$
In fact, since  $\dim\big(R/(\cap_{\dim F+1 \geq
d_{i+1}}I_F+\cap_{\dim F<d_i}I_F)\big)< d_i,$ and
$$D_i=\cap_{\dim
F+1\geq d_{i+1}}I_F/I_\Delta\simeq (\cap_{\dim F+1\geq
d_{i+1}}I_F+\cap_{\dim F< d_i}I_F)/\cap_{\dim F<d_i}I_F,$$
 we obtain
$$e(x_1,\dots,x_{d_i};D_i)=e(x_1,\ldots,x_{d_i};R/\cap_{\dim F<d_i}I_F).$$
Moreover, by the association law of multiplicity, we have
$$e(x_1,\ldots,x_{d_i};R/\cap_{\dim F<d_i}I_F)
=\sum_{\dim F=d_i-1}e(x_1,\ldots,x_{d_i};R/I_F).$$
Since
$R/I_F$ is a regular ring, $e(x_1,\ldots,x_{d_i};R/I_F)=\deg(x_1)\ldots \deg(x_{d_i})$; thus
$e(x_1,\ldots,x_{d_i};D_i) =\lambda_i \deg(x_1)\ldots
\deg(x_{d_i}).$
\end{proof}
\end{proposition}

We consider following examples, in which the number of vertices is
small. For the computation of these examples, we use CoCoA
\cite{cocoa} and Macaulay \cite{mac}.
\begin{example} For $n\leq 3$, $k[\Delta]$ are always sequentially Cohen-Macaulay
rings.
\begin{proof}
There are five cases (up to an isomorphism of $k[\Delta]$) as follows. \smallskip\\
i) $n\leq 2$ is trivial.\smallskip\\
ii) $k[\Delta]=R/(X_1X_2X_3)$:  It is easy to see that  $\un x
=(x_1=X_1+X_2,x_2=X_1+X_3)$ is a good system of parameters of
$k[\Delta]$. By a simple computation we have
$\ell(k[\Delta]/(x_1^2,x_2^2)k[\Delta])=12=2^2.3$; so $k[\Delta]$
is Cohen-Macaulay, since $\Delta$ has three facets  of
dimension $1$.
\smallskip\\
iii) $k[\Delta]=R/(X_2X_3)$: Similarly to (ii) we can find a good
system of parameters $\un x =(x_1=X_1,x_2=X_2+X_3)$. Then
$\ell(k[\Delta]/(x_1^2,x_2^2)k[\Delta])=8=2^2.2$. Since $\Delta$
has two facets  of dimension $1$, $k[\Delta]$ is
Cohen-Macaulay.\smallskip\\
iv) $k[\Delta]=R/(X_1X_2,X_2X_3,X_3X_1)$ is Cohen-Macaulay, since
$\dim k[\Delta]=\text{depth } k[\Delta]=1$.
\smallskip\\
v) $k[\Delta]=R/(X_1X_3,X_2X_3)$: The simplicial complex $\Delta$
has a facet of dimension $1$ and a facet of dimension $0$. There
is a good system of parameters $\un x = (x_1=X_2+X_3,x_2=X_1)$ of $k[\Delta]$ and
we have $\ell(k[\Delta]/(x_1^2,x_2^2)k[\Delta])=6=2.2+2$.
Therefore $k[\Delta]$ is sequentially Cohen-Macaulay but not
Cohen-Macaulay.
\end{proof}
\end{example}
\begin{example} For $n=4$, there is only one simplicial complex $\Delta$
 for which the ring $k[\Delta]$ is not
sequentially Cohen-Macaulay: Consider the simplicial complex
$\Delta$ given by two facets  $\{v_1,v_4\}, \{v_2,v_3\}$ of
dimension $1$. Then $k[\Delta]=R/(X_1X_2,X_1X_3,X_2X_4,X_3X_4).$
There is a good system of parameters $\un x=(x_1=X_1+X_2,
x_2=X_3+X_4)$ of $k[\Delta]$. We obtain by a simple computation
that $\ell(k[\Delta]/(x_1^2,x_2^2)k[\Delta])=9=2^2.2+1.$ So
$k[\Delta]$ is not sequentially Cohen-Macaulay. Moreover, it holds
$\ell(k[\Delta]/(x_1^{n_1},x_2^{n_2})k[\Delta])=2n_1n_2+1,$ and
therefore $k[\Delta]$ is a Buchsbaum ring. For all other
simplicial complexes $\Delta$, the ring $k[\Delta]$ are
sequentially Cohen-Macaulay. Par example, $\Delta$ is given by the
facets $\{X_1,X_2,X_3\},\{X_1,X_4\},$ $\{X_2,X_4\},\{X_3,X_4\}$.
Then $k[\Delta]=R/(X_1X_2X_4,X_1X_3X_4,X_2X_3X_4)$ has a good
system of parameters $\un
x=(x_1=X_3+X_4,x_2=X_1+X_2+X_3,x_3=X_1X_2)$. We have
$\ell(k[\Delta]/(x_1^2,x_2^2,x_3^2)k[\Delta])
=28=2^3.1.2+2^2.3,$ since $\Delta$ has one facet of dimension $2$
and three facets of dimension $1$. Therefore  $k[\Delta]$ is
sequentially Cohen-Macaulay.
\end{example}
\begin{example} For higher $n$, we can find many "bad" rings.
Consider a simplicial complex $\Delta$ over five vertices given by
the facets $\{v_1,v_2,v_3\},\{v_1,v_4,v_5\}$. Then
$$k[\Delta]=R/(X_2X_4,X_2X_5,X_3X_4,X_3X_5)=R/(X_2,X_3)\cap (X_4,X_5).$$
The ring $k[\Delta]$ has a good system of parameters $\un x=
(x_1=X_1;\ x_2=X_2+X_4;\ x_3=X_3+X_5)$. Moreover,
$$\ell(k[\Delta]/(x_1^{n_1},x_2^{n_2},x_3^{n_3})k[\Delta])=2n_1n_2n_3+n_1.$$
Since $\Delta$ has two facets  of dimension $2$, $k[\Delta]$ is
not sequentially Cohen-Macaulay. Note here that, although
$k[\Delta]/x_3k[\Delta]$ is a sequentially Cohen-Macaulay ring of
dimension $2$ and $x_3$ is a regular element of $k[\Delta]$, it
does not imply the sequentially Cohen-Macaulayness of $k[\Delta]$.
This is one of different features between sequentially
Cohen-Macaulayness and Cohen-Macaulayness.
\end{example}

\end{document}